\documentclass[a4paper]{article}

\usepackage[T1]{fontenc}
\usepackage[latin1]{inputenc}
\usepackage{latexsym}       
\usepackage{amsmath}
\usepackage{graphicx}
\usepackage{pslatex}
\usepackage{color}
\usepackage{epsfig}

\newtheorem{thm}{Theorem}[section]
\newtheorem{lem}[thm]{Lemma}
\newtheorem{cor}[thm]{Corollary}
\newtheorem{prop}[thm]{Proposition}
\newtheorem{conj}[thm]{Conjecture}

\newenvironment{monabstract}
  {\noindent\small\textbf{Abstract}.---\it}
  {\mbox{}\vspace{5mm}}

\newcommand{\metarem}[1]{}

\newcommand{\R}{\mathbf{R}}    
\newcommand{\Z}{\mathbf{Z}}    
\newcommand{\C}{\mathcal{C}}   
\newcommand{\tmax}{\tau_{\mathrm{max}}}   
\newcommand{\pt}{p_\tau}   
\newcommand{\conv}{\mathrm{Conv}} 
\newcommand{\extr}{\mathrm{Extr}} 
\newcommand{\len}{\mathrm{Len}} 
\newcommand{\T}{\mathcal{T}}   
\newcommand{\E}{\mathcal{E}}   
\newcommand{\M}{\mathcal{M}}   

\title{Triangulations of nearly convex polygons
       \footnote{\textsl{Math. Classif.}~:
          05B45, 52A10, 52A37, 52B55.}}
\author{Roland Bacher, Fr\'ed\'eric Mouton}

\begin{document}

\maketitle

\begin{monabstract}
Counting Euclidean triangulations with vertices in a finite set $\C$ of the convex hull 
$\conv(\C)$ of $\C$ is difficult in general, 
both algorithmically and theoretically.
The aim of this paper is to describe nearly convex polygons, a class of configurations for which
this problem can be solved to some extent.
Loosely speaking, a nearly convex polygon is an infinitesimal
perturbation of a weakly convex polygon (a convex polygon with edges subdivided
by additional points). 
Our main result shows that the
triangulation polynomial, enumerating all triangulations of a nearly convex polygon,
is defined in a straightforward way in terms of
polynomials associated to the ``perturbed'' edges.
\end{monabstract}
\section{Introduction}

Given a finite subset $\C$ of the Euclidean plane $\R^2$,
calculating the number of triangulations of the convex hull 
$\conv(\C)$  using only Euclidean triangles
with vertices in $\C$ seems to be difficult and has attracted
some interest, both from an algorithmic and a theoretical point of view,
see for instance \cite{A1}, \cite{A2}, \cite{AHN}, \cite{AK}, \cite{AAK},
\cite{An}, \cite{KZ}, \cite{S}, \cite{SS}.

An important and well understood special case is given by
the $n$ vertices of a strictly convex 
polygon. The associated number of triangulations is the Catalan number $C_{n-2}$.

In a first part of the paper, we consider convex polygons 
having perhaps collinear vertices, called weakly convex polygons.
We are not aware of the existence of formulae giving the number of triangulations for such polygons.
Thinking of weakly convex polygons as
strictly convex polygons with edges subdivided by additional vertices,
we call \textit{edges} the edges of the underlying strictly convex polygon.
Edges are thus maximal straight segments contained in the boundary of such polygons.

Denoting by $\C$ the set of all vertices, we define the 
\textit{weight} of an edge $E$ as the number of connected components of
$E\setminus (E\cap \C)$. Weights of successive edges form a finite sequence $a_1$, $a_2$... $a_l$ of
total sum $n=\# \C$. 

We show (Theorem~\ref{weaklyconvexmaxcase}) that there exists a sequence
of polynomials $p_m(t)$, $m\geq 1$, called \textit{maximal edge polynomials},
such that the number of maximal triangulations (\textit{i.e.} involving 
all vertices of $\C$) equals
\begin{equation}\label{formweaklyconvex}
  \tmax(\C)= \sum_{k\geq 2} b_k C_{k-2},
\end{equation}
where the coefficients $b_k$ are defined by $\prod_{i=1}^l p_{a_i} (t) = \sum_k b_k t^k$.

We deduce that the triangulation polynomial of a configuration (which takes into account
non-maximal triangulations) verifies formally the same formula as the previous one,
replacing maximal edge polynomials by \textit{complete edge polynomials}.
This has the perhaps surprising consequence that enumerative properties 
of triangulations do not depend on the particular cyclic order of the edges.

In a second part, we define \textit{nearly convex polygons} as small perturbations 
of weakly convex ones. Our main result, Theorem~\ref{nearlyconvexcase}, establishes the
existence of near-edge polynomials such that the 
previous formulae continue to hold.
Factorization of near-edges, a useful arithmetical property, allows 
classification of small near-edges.

Near-edge polynomials are difficult to compute in general except in a few special cases.
We plan to describe algorithms in a future paper dealing with computation aspects.
Some details are given in section~\ref{raq}.

This article is organized as follows.
We introduce first definitions and recall the strictly convex case in section~\ref{stopc},
prove formulae for weakly convex polygons in section~\ref{wcp}, expose the nearly convex setting
and prove our main result in section~\ref{ncp}. Finally, section~\ref{raq} contains a few
remarks and open problems.

\section{Triangulations of planar configurations}\label{stopc}

A {\it (planar) configuration of points} 
is a finite subset $\C=\{P_1,\ldots,P_n\}$ of the oriented plane $\R^2$.
We denote by $\conv(\C)$ the convex hull of $\C$ and by $\extr(\C)$ the set
of all extremal elements in $\C$ (a point $P\in\C$ is {\it extremal} if 
$\conv(\C\setminus\{P\})\neq \conv(\C)$).
The configuration $\C$ is said to be \textit{strictly convex} if
$\extr(\C)=\C$.
More generally, $\extr(\C)$ is the set of vertices of the strictly convex 
polygon formed by the convex hull of $\C$. 

A {\it triangulation} of a configuration $\C$ is a triangulation of 
its convex hull with vertices in $\C$, \textit{i.e.} a
finite set $\T=\{\Delta_1,\dots,\Delta_q\}$
of Euclidean triangles with vertices in $\C$ such that
$\conv(\C)=\cup_{i=1}^q\Delta_i$ and non-trivial intersections 
$\Delta_i\cap \Delta_j$ consist of a common
vertex or a common edge. A triangulation of $\C$ is  
\textit{maximal} if it involves all vertices of $\C$ (\textit{i.e.} each
point of $\C$ is a vertex of at least one triangle).
The number of maximal triangulations of $\C$
is denoted by  $\tmax(\C)$. If $\C$ is strictly convex,
it is well known that $\tmax(\C)$ is a \textit{Catalan number}:

\begin{thm} \label{strictlyconvexcase} All triangulations of a strictly convex $n-$gon
are maximal and their number is $C_{n-2}={2(n-2)\choose n-2}/(n-1)$.
\end{thm}

\noindent{\bf Sketch of proof:} 
A deformation argument shows that combinatorial properties of triangulations 
of a strictly convex $n-$gon depend only on $n$. We denote by $\tau_n$ the number of triangulations
of such a convex $n-$gon $P$.   
The choice of a marked edge $E$ in $P$ selects in every triangulation a unique triangle  
$\Delta$ containing $E$.
The two remaining edges of $\Delta$ determine two triangulated convex
polygons having respectively $k$ and $n+1-k$ edges for some integer $k$ such that $2\leq k \leq n-1$. 
This decomposition amounts to the recurrence relation $$\tau_n=\sum_{k=2}^{n-1}\tau_k\tau_{n+1-k} $$
holding for $n\geq 3$, using the convention $\tau_2=1$.
Therefore, the generating function $\sum_{n=2} \tau_n x^n$ satisfies a quadratic equation.
The classical resolution gives a formula, whose development in power series
yields the result. 
\hfill$\Box$

Triangulations of a general configuration $\C$ are not necessarily maximal and
enumerative properties are encoded by the triangulation polynomial
$\pt(\C) = \sum \tau_k(\C)s^k$, where  $\tau_k(\C)$ counts the number of triangulations
using exactly $k$ points. The polynomial $\pt(\C)$ has degree $n=\#\C$, with leading coefficient  $\tau_n(\C)=\tmax(\C)$
counting the number of maximal triangulations.
Its monomial of lowest degree $m=\#\extr(\C)$ corresponds to the $C_{m-2}$ triangulations of the convex hull
$\conv(\C)$ involving only extremal vertices.
Remark also that the average number of points of a triangulation is given by the logarithmic
derivative $f'(1)/f(1)$ of the triangulation polynomial $f(s)=\pt(\C)$.


Two configurations are \textit{isotopic} if they are related by 
a continuous deformation which preserves collinearity
and non-collinearity of triplets. Isotopic configurations have the same triangulation polynomial.

\section{Weakly convex polygons}\label{wcp}

\subsection{Definition and notations}

A configuration $\C$ is \textit{weakly convex} 
if it is contained in the boundary $\partial\conv(\C)$ of its convex hull.
We also call $\C$ a \textit{weakly convex polygon}.
We are not aware of a published formula giving the number of triangulations of such polygons. 

\begin{figure}[h]
\epsfxsize=\textwidth
\epsfbox{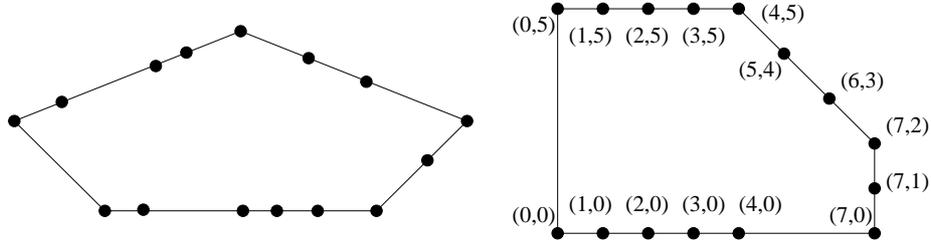}
\caption{Two weakly convex polygons with edge-weights $1,5,2,3,4$.}
\label{weighted-pentagon15234}
\end{figure}  

Weakly convex polygons can be seen as strictly convex polygons with additional vertices subdividing their edges.
We thus call \textit{edges} the segments joining two consecutive
extremal vertices of the underlying strictly convex polygon. 
An edge has \textit{weight} $a$ if it involves $a+1$ points of $\C$.
The weights of consecutive edges, in counterclockwise order, 
define, up to cyclic permutations, a finite sequence $a_1$, $a_2$... $a_l$ of
total sum $n=\# \C$ (see Figure~\ref{weighted-pentagon15234}). 
This sequence characterizes $\C$ up to isotopy. Thus all combinatorial properties of triangulations 
depend only on the sequence $a_1$, $a_2$... $a_l$ (up to cyclic permutations).
We denote by $\tau_k(a_1,a_2,\ldots,a_l)$ the number of corresponding triangulations
using exactly $k$ points of $\C$. This number is non-zero only for
$l\leq k\leq n$. We denote the triangulation polynomial of $\C$ by 
$\pt(a_1,a_2,\ldots,a_l)$.


The following notations will be useful.
The number $\tau_n(a_1,a_2,\ldots,a_l)$ of maximal triangulations is also denoted
by $\tau_{\mathrm{max}}(a_1,a_2,\ldots,a_l)$. 
We denote by  $P(a_1,a_2,\ldots,a_l)$ an arbitrary 
weakly convex polygon with $l$ edges of successive weights $a_1$, $a_2$... $a_l$.
Moreover, we use an exponential notation for indicating
several consecutive edges of weight $1$: 
we denote for instance by $P(1^{3},5,2)=P(1,1,1,5,2)$ a decagon with $5$ edges: three consecutive
edges of weight $1$, followed by an edge of weight $5$ and
a final edge of weight $2$.
We use the same notation for the number of triangulations: $\tau_6(1^{3},5,2)$
is the number of triangulations of $P(1^{3},5,2)$ involving $6$ vertices.


\subsection{Inclusion-exclusion principle}

Our first aim is the determination of the number of maximal triangulations for weakly convex polygons.
This can be achieved by reducing the problem to the case of strictly convex polygons
where formulae are known.
Replacing the first edge of weight $a_1$ by edges of weight $1$
leads to the following proposition.


\begin{prop}\label{propinclexcl}
  Given an integer $l\geq 3$ and  $l$ strictly positive integers $a_1$, $a_2$... $a_l$,
  we have 
  $$ \tau_{\mathrm{max}}(a_1,a_2,\ldots,a_l) = \sum_{k=0}^{\lfloor a_1/2\rfloor}
     (-1)^k{a_1-k\choose k}\tau_{\mathrm{max}}(1^{a_1-k},a_2,\ldots,a_l).$$  
\end{prop}

\begin{figure}[h]
\epsfxsize=\textwidth
\epsfbox{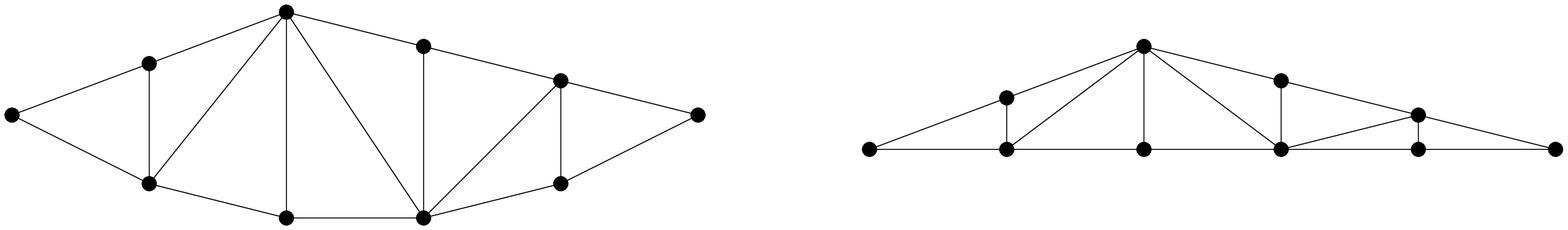}
\caption{Straightening a maximal triangulation without exterior triangles.}
\label{straightening}
\end{figure}

\noindent{\bf Proof.} We consider the set $\T$ of maximal triangulations of 
$P(1^{a_1},a_2,\ldots,a_3)$. A triangle of a triangulation is \textit{exterior} if
it involves two (necessarily adjacent) edges among the $a_1$ first edges (of weight $1$)
of $P(1^{a_1},a_2,\ldots,a_l)$. 
There is an obvious one-to-one correspondence between the set $\mathcal{R}$ of maximal triangulations
of $P(1^{a_1},a_2,\ldots,a_l)$ without exterior triangles 
and the set of maximal triangulations of $P(a_1,a_2,\ldots,a_l)$, by continuously straightening
the set formed by the first $a_1$ $1-$edges (see Figure~\ref{straightening}). It is thus sufficient to
enumerate $\mathcal{R}$. It is easier to enumerate the complementary set $\T\setminus\mathcal{R}$.


Denoting by $\E$ the set of all $a_1-1$ possible exterior triangles and, for $\Delta\in \E$,
$\T_\Delta$ the set of triangulations containing $\Delta$, the set  $\T\setminus\mathcal{R}$
is the union of the sets $\T_\Delta$. We enumerate this set by the inclusion-exclusion principle:
$$  \#(\T\setminus\mathcal{R}) = \sum_{k=1}^{\lfloor a_1/2\rfloor}
                         (-1)^{k-1}
                         \sum_{\{\Delta_1,\Delta_2,\ldots,\Delta_k\}}
                            \#(\T_{\Delta_1}\cap \T_{\Delta_2}\cap \cdots \cap  \T_{\Delta_k})   .$$
The upper bound ${\lfloor a_1/2\rfloor}$ in the summation is due to the fact that a triangulation
contains at most ${\lfloor a_1/2\rfloor}$ exterior triangles.

\begin{figure}[h]
\epsfxsize=\textwidth
\epsfbox{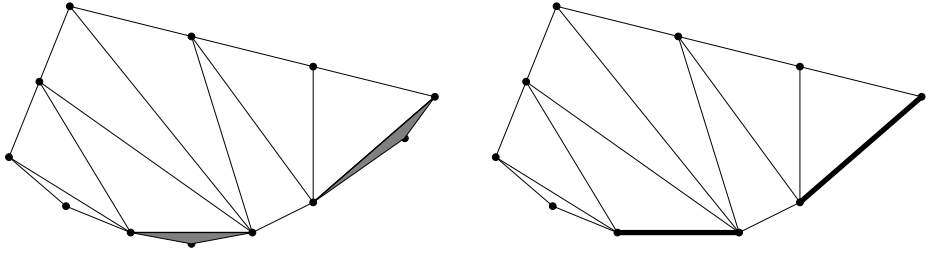}
\caption{Erasing some exterior triangles and keeping marks.}
\label{2decorated}
\end{figure}

It remains to enumerate the intersections. Fix $k\leq a_1/2$. If some triangles among 
$\Delta_1$, $\Delta_2$... $\Delta_k$ have non-disjoint interiors, the intersection is empty.
Otherwise, we associate to each element of $\T_{\Delta_1}\cap \T_{\Delta_2}\cap \cdots \cap  \T_{\Delta_k}$
a maximal triangulation of the polygon $P(1^{a_1-k},a_2,\ldots,a_l)$ by erasing triangles
$\Delta_1$, $\Delta_2$... $\Delta_k$. We also keep track of erased triangles
by marking the remaining edge for each of these triangles (see Figure~\ref{2decorated}).
This defines a map $\varphi_{\Delta_1, \Delta_2,\ldots,\Delta_k}$ from the set
$\T_{\Delta_1}\cap \T_{\Delta_2}\cap \cdots \cap  \T_{\Delta_k}$ to the set $\M_k$ of all
triangulations of $P(1^{a_1-k},a_2,\ldots,a_l)$ with $k$ marked edges among the first 
$a_1-k$ edges of weight $1$. 
This map is obviously injective: we can reconstruct
the initial triangulation by gluing triangles onto the marked edges.
Moreover, the union of images of $\varphi_{\Delta_1, \Delta_2,\ldots,\Delta_k}$, 
for all ``admissible'' $k-$tuples of exterior triangles is clearly
a disjoint union, by the same remark as for injectivity, and fills $\M_k$
for exactly the same reason: the reconstruction is unique and is always possible.
Thus,
$$  \sum_{\{\Delta_1,\Delta_2,\ldots,\Delta_k\}}
        \#(\T_{\Delta_1}\cap \T_{\Delta_2}\cap \cdots \cap  \T_{\Delta_k}) 
    = \#\M_k = {a_1-k \choose k} \tau_{\mathrm{max}}(1^{a_1-k},a_2,\ldots,a_l)  .$$
As $\#\T = \tau_{\mathrm{max}}(1^{a_1},a_2,\ldots,a_l)$, the proof is achieved.
\hfill$\Box$

\subsection{Maximal triangulations}\label{ssmaxtri}

Since combinatorial properties of triangulations of  weakly convex polygons 
are invariant under cyclic permutations of edges,
we can ``break up'' all edges by iterating Proposition~\ref{propinclexcl}. 
Hence, $\tau_{\mathrm{max}}(a_1,a_2,\ldots,a_l)$
can be successively written as
\begin{eqnarray*} 
     &  & \sum_{k_1=0}^{\lfloor a_1/2\rfloor}
           (-1)^{k_1}{a_1-k_1\choose k_1}
           \tau_{\mathrm{max}}(1^{a_1-k_1},a_2,\ldots,a_l)   \\
     & = &  \sum_{k_1=0}^{\lfloor a_1/2\rfloor} \sum_{k_2=0}^{\lfloor a_2/2\rfloor}
            (-1)^{k_1+k_2}{a_1-k_1\choose k_1}{a_2-k_2\choose k_2}
           \tau_{\mathrm{max}}(1^{(a_1-k_1)+(a_2-k_2)},a_3,\ldots,a_l)   \\
     & = & \ldots\\ 
     & = &  \sum_{k_1,k_2,\ldots,k_l}
               \left( (-1)^{\sum_{i=1}^l k_i}\prod_{i=1}^l{a_i-k_i\choose k_i} \right)
                \tau_{\mathrm{max}}(1^{\sum_{i=1}^l a_i-k_i}) \\   
     & = &  \sum_{k_1,k_2,\ldots,k_l}
               \left( (-1)^{\sum_{i=1}^l k_i}\prod_{i=1}^l{a_i-k_i\choose k_i} \right)
                 C_{(\sum_{i=1}^l a_i-k_i)-2}. \\   
\end{eqnarray*}
We thus obtain $\tau_{\mathrm{max}}(a_1,a_2,\ldots,a_l)$ as a linear combination 
$\sum_{j\geq 2} b_j C_{j-2}$ of Catalan numbers.
The coefficients $b_j$ are given by 
$$\sum_{j\geq 2} b_j t^j =
  \sum_{k_1=0}^{\lfloor a_1/2\rfloor} \cdots \sum_{k_l=0}^{\lfloor a_l/2\rfloor}
               \left( (-1)^{\sum_{i=1}^l k_i}\prod_{i=1}^l{a_i-k_i\choose k_i} \right)
                 t^{\sum_{i=1}^l a_i-k_i} =
  \prod_{i=1}^l p_{a_i}$$
with
\begin{equation}\label{maxedgepol}
   p_{m}= \sum_{k=0}^{\lfloor m/2\rfloor} (-1)^k{m-k\choose k}t^{m-k}
\end{equation}
defining
the sequence $(p_m)m\geq 1$ of \textit{maximal edge-polynomials}.

We consider the generating function 
\begin{equation}\label{gencatalan}
  G_C(t)=\sum_{k\geq 0} C_{k}t^k = \sum_{k\geq 0} {2k \choose k}\frac{t^k}{k+1}  
\end{equation}
for the sequence of Catalan numbers (corresponding to the analytic expression $\frac{1-\sqrt{1-4t}}{2t}$).
Given a polynomial $p(t)= \sum_{k\geq 2} \alpha_k t^k$, we define
\begin{equation}\label{scalprod}
  \big\langle p(t),t^2G_C(t)\big\rangle_t = \sum_{k\geq 2} \alpha_k C_{k-2}.
\end{equation}
This ``umbral'' notation is suggested by the fact that polynomials and formal
power series are mutually dual.

We have obtained the following result concerning the number of maximal triangulations:
\begin{thm}\label{weaklyconvexmaxcase}
  Given natural numbers $l\geq 3$ and $a_1$, $a_2$... $a_l\geq 1$, we have 
  $$ \tau_{\mathrm{max}}(a_1,a_2,\ldots,a_l) = \big\langle \prod_{i=1}^l p_{a_i}(t),t^2G_C(t) \big\rangle_t,$$  
  using the notations of 
  formulae~\ref{maxedgepol}, \ref{gencatalan} and~\ref{scalprod} above.
\end{thm}

\subsection{Examples and remarks}

The first few maximal edge-polynomials are
$$\begin{array}{ll}
\displaystyle p_1=t&
\displaystyle p_5=t^5-4t^4+3t^3\\
\displaystyle p_2=t^2-t&
\displaystyle p_6=t^6-5t^5+6t^4-t^3\\
\displaystyle p_3=t^3-2t^2&
\displaystyle p_7=t^7-6t^6+10t^5-4t^3\\
\displaystyle p_4=t^4-3t^3+t^2\qquad&
\displaystyle p_8=t^8-7t^7+15t^6-10t^5+t^4\\
\end{array}$$

\noindent{\bf Example.} The two weakly convex polygons $P(1,5,2,3,4)$
of Figure~\ref{weighted-pentagon15234} have 
$$\begin{array}{cl}
    \displaystyle   & \langle p_1\ p_5\ p_2\ p_3\ p_4,t^2G_C(t) \rangle_t\\ 
    \displaystyle = & \langle t(t^5-4t^4+3t^3)(t^2-t)(t^3-2t^2)(t^4-3t^3+t^2),t^2G_C(t)\rangle_t\\ 
    \displaystyle = & \langle t^{15}-10t^{14}+39t^{13}-75t^{12}+74t^{11}-35t^{10}+6t^9,t^2G_C(t)\rangle_t\\ 
    \displaystyle = & C_{13}-10C_{12}+39C_{11}-75C_{10}+74C_{9}-35C_8+6C_7\\ 
    \displaystyle = & 7429000-10\cdot 208012+39\cdot 58786-75\cdot 16796\\
    \displaystyle   & \quad +74\cdot 4863-35\cdot 1430+6\cdot 429\\ 
    \displaystyle = & 8046
\end{array}$$
maximal triangulations. 

\noindent{\bf Remark.} We have $\tau_{max}(n,m,1^2)={n+m\choose n}$ for all $n,m\geq 1$.\\
Indeed, triangles in a maximal triangulation of $P(1,n,1,m)$
are linearly ordered and in one-to-one correspondence
with the $n+m$ "segments" of the two opposite "long" edges.
Gluing a triangle onto one of the two edges of length one, we have the formula
$\tmax(n+1,m+1,1) = {n+m \choose n}$.

\noindent{\bf Remark.} The function $n\longmapsto f(k,n)= \tmax(1^{2+k},n)$ is polynomial of degree $k$.\\
First, a classification of the triangulations of $P(1^{2+k},n)$ according to the third vertex
of the last triangle based on the edge of length $n$ gives the formula
$$ f(k,n)= \sum_{l=0}^k C_{k-l} f(l,n-1)  .$$ 
The result follows then from an induction on $k$. It is obvious that $f(0,n)=1$ for all $n$.
Suppose that $f(l,n)$ is polynomial in $n$ of degree $l$ for every $l\leq k$.
Using $C_0=1$, the difference 
$$f(k+1,n)-f(k+1,n-1)= \sum_{l=0}^{k} C_{k+1-l} f(l,n-1)$$ 
is then a polynomial of degree $k$ and a sum over $n$ implies the result.

\noindent{\bf Remark.} The sequence of numbers of maximal triangulations of the weakly convex polygons 
$P(2,2,2)$, $P(2,2,2,2)$... starts as~:\\
\centerline{4, 30, 250, 2236, 20979, 203748, 2031054, 20662980, 213679114, 2239507936...}
(see sequence A86452 of~\cite{OEIS}).

\noindent{\bf Remark.}
  Maximal edge polynomials can also be defined recursively by $p_0=1$, $p_1=t$ and
  $p_m=t(p_{m-1}-p_{m-2})$ and are related to Fibonacci
  numbers (the closely related polynomials
  $\sum_{k}{m-k\choose k}x^k$ are also called ``Fibonacci polynomials'').

\metarem{Lien avec les Tchebichev}

\metarem{$\alpha_k=\lim_{n\rightarrow \infty} \tmax{k^n}^{\frac{1}{n}}$ existe et
vérifie $\alpha_{k+l}\leq \alpha_k \alpha_l$? Quel est sa croissance? Nombres algébriques?
$\alpha_2=4$}

\subsection{Non-maximal triangulations}

An arbitrary (\textit{i.e.} not necessarily maximal) triangulation of a weakly convex
polygon $P(a_1,\ldots,a_l)$ is a maximal triangulation of a subset
involving all extremal vertices of the weakly convex configuration
$P(a_1,\ldots,a_l)$.
It amounts thus to the choice, for every $1\leq i\leq l$, of a number $1\leq b_i\leq a_i$ and of $b_i-1$ points 
among the $a_i-1$ interior points of the $i-$th edge, followed by the choice of a triangulation of $P(b_1,\ldots,b_l)$.
  

The triangulation polynomial of $P(a_1,\ldots,a_l)$ is thus
given by
$$
\begin{array}{l}
     \displaystyle \sum_{1\leq b_i\leq a_i}
                \left(\prod_{j=1}^l {a_i-1\choose b_i-1}\right)
                \big\langle \prod_{j=1}^lp_{b_j}(t),t^2G_C(t) \big\rangle_t
                \ s^{\sum_{j=1}^l b_j}\\
     \displaystyle 
       \qquad = \big\langle
                  \prod_{j=1}^l\sum_{b_j=1}^{a_j}{a_j-1\choose b_j-1}p_{b_j}(t)s^{b_j},
                  t^2G_C(t)
                \big\rangle_t\\ 
     \displaystyle 
       \qquad = \big\langle\prod_{j=1}^l{\overline p}_{a_j},t^2G_C(t) \big\rangle_t,
\end{array}
$$
where the \textit{complete edge-polynomials} 
$\overline{p}_m\in\Z[s,t]$ are defined as
\begin{equation}\label{compedgepol}
  {\overline p}_m=\sum_{k=1}^m {m-1\choose k-1}p_k(t)s^k .
\end{equation}
We have proved

\begin{thm}\label{weaklyconvexpolcase}
  Given natural numbers $l\geq 3$ and $a_1$, $a_2$... $a_l\geq 1$, 
  the triangulation polynomial of $P(a_1,a_2,\ldots,a_l)$ is
  $$ \pt(a_1,a_2,\ldots,a_l) = \sum_{k}\tau_k(a_1,a_2,\ldots,a_l)s^k = 
      \big\langle \prod_{i=1}^l \overline{p}_{a_i}(t),t^2G_C(t) \big\rangle_t,$$  
  using the notations of 
  formulae~\ref{gencatalan}, \ref{scalprod} and~\ref{compedgepol} above.
\end{thm}

An immediate consequence is the following slightly surprising fact:

\begin{cor}
  Enumerative properties of triangulations
  for weakly convex polygons 
  do not depend on the particular cyclic order of edge weights.
\end{cor}

The first few complete edge-polynomials are
$$\begin{array}{l}
\displaystyle {\overline p}_1=p_1\ s=t\ s\ ,\\
\displaystyle {\overline p}_2=p_2\ s^2+p_1\ s=(t^2-t)\ s^2+t\ s\ ,\\
\displaystyle {\overline p}_3=p_3\ s^3+2p_2\ s^2+p_1\ s=(t^3-2t^2)\ s^3+
2(t^2-t)\ s^2+t\ s\ ,\\
\displaystyle {\overline p}_4=p_4\ s^4+3p_3\ s^3+3p_2\ s^2+p_1\ s\\
\displaystyle \qquad =(t^4-3t^3+t^2)\ s^4+3(t^3-2t^2)\ s^3+
3(t^2-t)\ s^2+t\ s\ .\end{array}$$

\noindent{\bf Example.} The triangulation polynomial 
$\pt(1,5,2,3,4)$ of the two weakly convex 
polygons of Figure \ref{weighted-pentagon15234} equals
$$\begin{array}{cl}
     \displaystyle   & \langle 
                       {\overline p}_1\ {\overline p}_5\ {\overline p}_2\ {\overline p}_3\ {\overline p}_4,t^2G_C(t)
                     \rangle_t\\ 
     \displaystyle = & 8046\,{s}^{15}+37250\,{s}^{14}+77467\,{s}^{13}+95364\,{s}^{12}+77048\,{s}^{11}\\
     \displaystyle   & \quad +42776\,{s}^{10}+16584\,{s}^{9}+4460\,{s}^{8}+805\,{s}^{7}+90\,{s}^{6}+5\,{s}^{5}\ .
\end{array}$$

\section{Nearly convex polygons}\label{ncp}

Nearly convex polygons are small perturbations of weakly convex polygons and form
the correct framework for generalizing 
Theorem~\ref{weaklyconvexmaxcase} and Theorem~\ref{weaklyconvexpolcase}
We give first a definition of near-edges, which are small deformations of edges, and introduce
nearly convex polygons and, using the formalism of ``roofs'', the associated near-edge polynomials. 
Then we state and prove the main theorem. A last subsection describes 
factorization properties of near-edges and gives a classification of small near-edges.

\subsection{Near-edges}

In order to describe deformations of an edge $E$ in a weakly convex polygon $P=P(n,\dots)$,
we choose coordinates such that $P$ is contained in the upper half-plane
$y\geq 0$ and $E$ is a subset of the boundary $y=0$. 
This leads to the following definition.

A \textit{near-edge} of \textit{weight} $n$ (or an $n-${\it near-edge}) is a sequence $E$ of 
$n+1$ points
\begin{equation}\label{nearedgecoord}
  (P_0,P_1,\ldots,P_n)=
    \left(
      \left(\begin{array}{c}x_0\\y_0\end{array}\right),
      \left(\begin{array}{c}x_1\\y_1\end{array}\right),
      \ldots,
      \left(\begin{array}{c}x_{n}\\y_n\end{array}\right)
    \right)
    \in\left({\mathbf R}^2\right)^{n+1}
\end{equation}
such that $x_0<x_1<\dots<x_{n}$ and $y_0=y_{n}=0$. We 
consider the complete order $P_0<P_1<P_2<\dots<P_n$ on $E$ and 
call $P_0$, respectively $P_n$, the {\it initial}, respectively {\it final},
vertex of the near-edge $E$. We denote a
near-edge $E$ either by the sequence $E=(P_0,\dots,P_n)$ 
of its $n+1$ points or by the  real matrix 
$$E=\left(\begin{array}{ccccccc}x_0&x_1&\dots&x_{n-1}&x_n\\ 0&y_1&\dots&y_{n-1}&0\end{array}\right)$$
of size $2\times (n+1)$.
 
A continuous deformation of near-edges, 
which preserves collinearity and non-collinearity of all triplets of points,
is called an {\it isotopy}. Two near-edges
joined by an isotopy are {\it isotopic}.

Given a near-edge $E$ with points $P_i\in\R^2$ as above,
we denote by $E^\epsilon$ the near-edge with points
$$P_k(\epsilon)=\left(\begin{array}{c}x_k\\\epsilon
\ y_k\end{array}\right)=\left(\begin{array}{cc}1&0\\0&\epsilon\end{array}
\right)P_k ,$$
for $0\leq k \leq n$.

In particular $E=E^1$ and all near-edges $E^\epsilon$ are isotopic 
for $\epsilon>0$.

\subsection{Nearly convex polygons}

Let $P$ be a strictly convex polygon with $l$ extremal vertices
$V_0$, $V_1$, $V_2$\ldots $V_l=V_0$,
appearing in counterclockwise order around the boundary $\partial P$ of $P$.

Given a sequence $E_1,\ldots,E_l$
where $E_i=(P_{i,0},P_{i,1},\ldots,P_{i,n_i})$ 
is an $n_i-$near-edge, we denote by $G(E_1,\ldots,E_l|P)$
the unique configuration obtained by gluing the 
$n_i-$near-edge $E_i$, rescaled suitably by an 
orientation-preserving similitude, along the oriented
edge of $P$ which starts at $V_{i-1}$ and ends at $V_i$. 

More precisely, the gluing map $\varphi_i$ is the unique orientation-preserving 
similitude of ${\mathbf R}^2$ such that 
$$\varphi_i(P_{i,0})=V_{i-1}\hbox{ and }\varphi_i(P_{i,n_i})=V_{i}\ .$$
The configuration $G(E_1,\ldots,E_l|P)$ is the set of points
$\cup_{i=1}^l \varphi_i(E_i) \subset {\mathbf R}^2$.

We have the following result which we state without proof.

\begin{prop}\label{isotopyclass} 
  Consider a strictly convex $l-$gon $P$ and $l$ near-edges $E_1,\ldots,E_l$,
  \begin{enumerate}
     \item The configurations $G(E_1^{\epsilon_1},\ldots,E_l^{\epsilon_l}|P)$ are all isotopic
           for all $\epsilon_i>0$ small enough. This defines an isotopy class associated
           to the near-edges and the polygon.
     \item Given a second strictly convex $l-$gon $Q$, the configurations $G(E_1^{\epsilon_1},\ldots,E_l^{\epsilon_l}|P)$
           and $G(E_1^{\epsilon_1},\ldots,E_l^{\epsilon_l}|Q)$ are isotopic for all $\epsilon_i>0$ small enough.
           Hence, the isotopy class above does not depend on $P$.
     \item The isotopy class defined in this way depends only on the isotopy classes of the near-edges.
  \end{enumerate}
\end{prop}

A \textit{nearly convex polygon} is 
a configuration in the isotopy class associated by Proposition~\ref{isotopyclass} 
to a sequence $E_1,\ldots,E_l$ of near-edges. 
For the sake of convenience, the isotopy class itself is also called a \textit{nearly convex polygon}.
As far as combinatorial properties of triangulations are concerned, all the configurations of the class
are equivalent and we denote by $P(E_1,\ldots,E_l)$ any such configuration. 
This notation is a natural extension of the notation already used for weakly convex polygons,
with integers $n$ representing weighted edges of weakly convex polygons. 

\begin{figure}[h]
\epsfxsize=\textwidth
\epsfbox{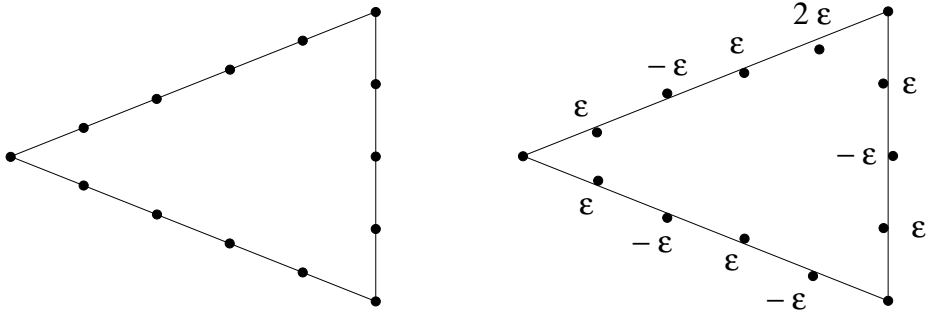}
\caption{Perturbating $P(5,4,5)$ into a nearly convex polygon.}
\label{exple1}
\end{figure}  
 
One can think of a nearly convex polygon as a small perturbation of the
configuration associated to a weakly convex polygon. Figure~\ref{exple1}
shows the weakly convex polygon $P(5,4,5)$ and
a nearly convex polygon obtained 
by moving slightly non-extremal vertices perpendicularly to the
three corresponding edges (with a hopefully 
evident notation indicating the perturbation). 
This nearly convex polygon is isotopic to $P(E_a,E_b,E_c)$ where
$$\begin{array}{l}
\displaystyle
E_a=\left(\begin{array}{rrrrrr}0&1&2&3&4&5\\0&1&-1&1&-1&0\end{array}\right)\\
\displaystyle 
E_b=\left(\begin{array}{rrrrr}0&1&2&3&4\\0&1&-1&1&0\end{array}\right)\\
\displaystyle
E_c=\left(\begin{array}{rrrrrr}0&1&2&3&4&5\\0&2&1&-1&1&0\end{array}\right).
\end{array}$$

\metarem{La representation intégrale a-t-elle un intérêt?}
%

\subsection{Edge-type triangles}

When considering a triangulation $\T$ of a nearly convex polygon $P$,
and a near-edge $E$ of $P$,
some triangles of $\T$ have their three vertices in $E$
(we identify here the near-edge and its realization in $P$).
We call such a triangle $\Delta$ an \textit{edge-type triangle} of $E$ and we write $\Delta\subset E$.
The following technical lemma is the key of our main result.

\begin{lem}\label{keylemma}
  Given a triangulation $\T$ of a nearly convex polygon $P$
  and a near-edge $E$ of $P$, each point of $E$ that is a vertex of $\T$ 
  belongs to $\partial \conv(P)$ or to at least one edge-type triangle of $E$.  
\end{lem}

\noindent{\bf Proof.} 
Recall that nearly convex polygons are ``small'' perturbations of 
weakly convex polygons. 
For the fixed particular realization of the near-edge $E$ in the current polygon,
denote by $E^\epsilon$ the realization obtained
by setting all points of $E$ closer to the edge $[P_0,P_n]$
by multiplying the distances to this edge by $\epsilon\in (0,1)$.
We can then replace $E$ by $E^\epsilon$ ($\epsilon\in (0,1)$)
whithout changing the structure of the triangulation (Proposition~\ref{isotopyclass}).

Let us prove the lemma by contradiction: suppose that there exists a point $M$ of $E$ which is an interior 
vertex of $\T$ and which belongs to no edge-type triangle of $E$.
Hence, there are at least three edges of the triangulation starting at $M$ and forming adjacent angles less than $\pi$.
If no such edge joins two vertices of $E$,
a substitution of $E$ by $E^\epsilon$, for appropriately small $\epsilon$,  
moves $M$ so close to the segment $[P_0,P_n]$ that one of these angles
becomes greater than $\pi$.
Therefore, one at least of these edges has its other endpoint $N\in E$. 
Suppose for convenience, but without loss of generality, that $M<N$.
Then consider the next segment of the triangulation starting from $M$,
according to the clockwise order around $M$, and denote by $Q$ its other endpoint. 
The asumptions on $M$ imply that $Q\not\in E$.
Replacing now $E$ by $E^\epsilon$, for $\epsilon$ sufficently small,
gives an angle $\widehat{QMN} > \pi$ for the triangle $(QMN)$, which is a contradiction.
\hfill$\Box$

\subsection{Roofs}

Let us define the set
$$R(\T,E)=\partial \left(\conv(P)\setminus \bigcup_{\Delta\subset E}\Delta\right)\cap \conv(E),$$
which will be proved to be a piecewise linear path 
separating the edge-type triangles contained in $E$ from the remaining triangles of $\T$.

The idea for counting triangulations of a nearly convex polygon 
is to classify triangulations according to the paths thus obtained
from all near-edges and to enumerate all triangulations giving rise to such
a set of paths.
The following definition is useful for the description of all possibilities.

A \textit{partial roof} with \textit{length} $\len(R)=k$ of $E$ is a 
piecewise linear path $R$ starting at the initial vertex and ending at the final vertex of $E$,
whose $k+1$ vertices are elements of $E$ in increasing order:
$P_{j_0}=P_0,P_{j_1},\dots,P_{j_{k-1}},P_{j_k}=P_n$.
This partial roof is denoted by
$$R=[P_{j_0}=P_0,P_{j_1},\dots,P_{j_{k-1}},P_{j_k}=P_n]$$ 

\begin{lem}\label{lempartroof}
  Given a triangulation $\T$ of a nearly convex polygon $P$
  and a near-edge $E$ of $P$, the set
  $R(\T,E)$ is a partial roof of $E$.
\end{lem}

\noindent{\bf Proof.} 
Let us introduce the set $P'=\overline{\conv(P)\setminus \bigcup_{\Delta\subset E}\Delta}=\bigcup_{\Delta\not\subset E}\Delta$.
Since $\conv(P)\setminus\conv(E)$ is a connected subset of $P'$ which intersects all triangles $\Delta\not\subset E$,
$P'$ is connected. Hence, it is a (generally non convex) polygon with vertices in $\T$.

Considering the inclusions $\conv(P)\setminus\conv(E)\subset P'\subset \conv(P)$, the boundary of $P'$
coincides with the boundary of $\conv(P)$ outside of $\conv(E)$. Hence, $R(\T,E)$ is a piecewise linear path
with vertices $P_{j_0},P_{j_1},\dots,P_{j_{k-1}},P_{j_k}$ in $E$,
whose orientation can be chosen such that $P_{j_0}=P_0$ and $P_{j_k}=P_n$ (where $P_0$ and $P_n$ are
the initial and final vertices of $E$).

Let us now prove by contradiction that this sequence of points of $E$ is increasing. Otherwise,
consider the first decreasing step: $P_{j_{i+1}}<P_{j_{i}}$. 
There are two cases.

Suppose first that the point $P_{j_{i+1}}$ is above the line $L$ defined by $P_{j_{i-1}}$ and $P_{j_{i}}$, according to
the standard coordinates of the near-edge $E$.
By definition
of $R(\T,E)$, the point $P_{j_i}$ is linked by an edge of $\T$ to a point $Q\in P\setminus E$ which crosses
the segment $[P_{j_{i-1}},P_{j_{i+1}}]$. As in the proof of Lemma~\ref{keylemma}, a substitution
of $E$ by $E^\epsilon$, for $\epsilon$ sufficiently small, removes this crossing, which is in
contradiction with Proposition~\ref{isotopyclass}.

Suppose now that the point $P_{j_{i+1}}$ is below the line $L$. The point
$P_{j_{i+1}}$ is linked by an edge of $\T$ to a point $Q'\in P\setminus E$. A substitution
of $E$ by $E^\epsilon$, for $\epsilon$ sufficiently small, creates a crossing between this edge and the 
edge $[P_{j_{i-1}},P_{j_{i}}]$, which is in
contradiction with Proposition~\ref{isotopyclass}.
\hfill$\Box$

The definition of a partial roof refers to 
the coordinate representation of $E$: the ``sky'' (which corresponds to the interior of the
nearly convex polygon)
is ``above'' $E$. We are interested in points and triangles ``sheltered'' by a
partial roof.

In order to define properly the region sheltered by a partial roof, we consider again
coordinates (formula~\ref{nearedgecoord}) of the near-edge $E$:
we define the \textit{lower boundary} $\partial^- E$ of $E$ as the  piecewise-linear path
$\partial^- E=(\partial\conv(E))\cap\{(x,y)\ \vert\ y\leq 0\}$, which is the ``lowest'' 
possible partial roof.
Remark that each partial roof $R$, and in particular the lower boundary, is the graph 
of a piecewise-affine function $f_R: [x_0,x_n]\longrightarrow \R$.
The graph $\partial^- E$ is below each other partial roof: we have
$f_{\partial^- E}(x)\leq f_R(x)$, for all $x$ in $[x_0,x_n]$. 
The region \textit{sheltered} by the partial roof $R$ is then the (generally non convex) subset $S(E,R)$ enclosed by    
$\partial^- E$ and $R$~: 
$$S(E,R)=\{ (x,y)| x\in [x_0,x_n], f_{\partial^- E}(x)\leq y \leq f_R(x)\}.$$
We use the same notations $\partial^- E$ and $S(E,R)$ for denoting the corresponding subsets
of a realization of $E$ in a nearly convex polygon $P$.

\begin{lem}\label{lemtriangunion}
  Given a triangulation $\T$ of a nearly convex polygon $P$
  and a near-edge $E$ of $P$, $\T$ induces a triangulation on the sheltered region:
  $$S(E,R(\T,E))=\bigcup_{\Delta \subset E} \Delta.$$  
\end{lem}

\noindent{\bf Proof.} It is a consequence of the proof of Lemma~\ref{lempartroof}.
\hfill$\Box$

A partial roof $R$ is a \textit{roof} if $E\subset S(E,R)$. 
If $\T$ is a maximal triangulation of $P$, $R=R(\T,E)$ is a roof by Lemma~\ref{keylemma}
and $\T$ induces a triangulation on $S(E,R)$ whose vertices are exactly the points of $E$ by Lemma~\ref{lemtriangunion}.
Enumerating induced triangulations on $S(E,R)$ leads to the definition of
maximal near-edges polynomials. 

\subsection{Near-edge polynomials}

We define the \textit{maximal polynomial} of the near-edge $E$ by
\begin{equation}\label{maxnearedgepol}
  p_E=\sum_{R\mathrm{\ roof\ of\ }E} \tmax(E,R)\ p_{\len(R)}\ \ \in\Z[t],
\end{equation} 
where $\tmax(E,R)$ denotes the number of triangulations of $S(E,R)$
involving exactly all points of $E\cap S(E,R)$ and 
the polynomials $p_m$ are the maximal edge polynomials given by formula~\ref{maxedgepol} of subsection~\ref{ssmaxtri}.

Starting from a triangulation $\T$ which is not maximal, the path $R$ induced on a near-edge $E$
is in general only a partial roof and $\T$ induces not necessarily a maximal triangulation
of the sheltered region $S(E,R)$.
This suggest to introduce ``sub-near-edges'' in order
to deal with these difficulties.

Let $V^-(E)=\extr(E)\cap\partial^- E$ denote the set of extremal points of the 
lower boundary. A {\it $k-$sub-near-edge} of $E$ is an increasing subsequence 
$E'\subset E$ of $k+1\leq n+1$ elements containing the set 
$V^-(E)$. In particular, any sub-near-edge $E'$ of $E$ has initial 
vertex $P_0$, final vertex $P_n$ and verifies $\partial^-E'=\partial^-E$. 
We note $E'<E$.

\begin{figure}[h]
\epsfxsize=\textwidth
\epsfbox{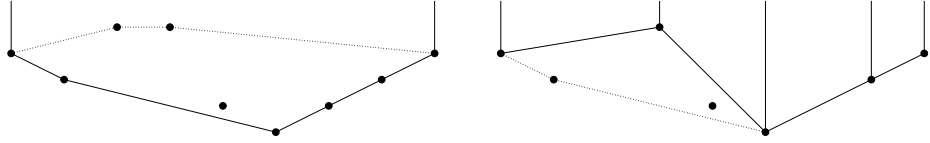}
\caption{An $8-$near-edge and a roof of length $4$ of a $6-$sub-near-edge.}
\label{8near-edge}
\end{figure}  

\noindent{\bf Example.} The left half of Figure~\ref{8near-edge}
displays the $8-$near-edge
$$E=(P_0,\dots,P_8)=\left(\begin{array}{rrrrrrrrr}
0&1&2&3&4&5&6&7&8\\ 
0&-1&1&1&-2&-3&-2&-1&0\end{array}\right)\ .$$
We have $V^-(E)=(P_0,P_1,P_5,P_8)$
and $E$ has $2^5$ sub-near-edges obtained by removing any subset of 
vertices among $\{P_2,P_3,P_4,P_6,P_7\}$
from $E$. The right half of Figure~\ref{8near-edge} 
shows the roof $R=[P_0,P_3,P_5,P_7,P_8]$
of the sub-near-edge $E'$ of $E$
defined by $(P_0,P_1,P_3,P_4,P_5,P_7,P_8)$.

The {\it complete polynomial} ${\overline p}(E)$ of an $n-$near-edge
$E=(P_0,\dots,P_n)$ is defined as
\begin{equation}\label{compnearedgepol}
  \overline{p}_E=\sum_{m=1}^n \sum_{\genfrac{}{}{0pt}{}{E'< E}{\#E'=m}} p_{E'}s^m.
\end{equation}

{\bf Example.} We compute the complete polynomial 
${\overline p}_{E_a}$ of the near-edge 
$$E_a=(P_0,\dots,P_5)=
\left(\begin{array}{rrrrrr}0&1&2&3&4&5\\0&1&-1&1&-1&0\end{array}\right)
$$
involved in the nearly convex polygon $P(E_a,E_b,E_c)$ of Figure~\ref{exple1}.

\begin{figure}[h]
\epsfxsize=\textwidth
\epsfbox{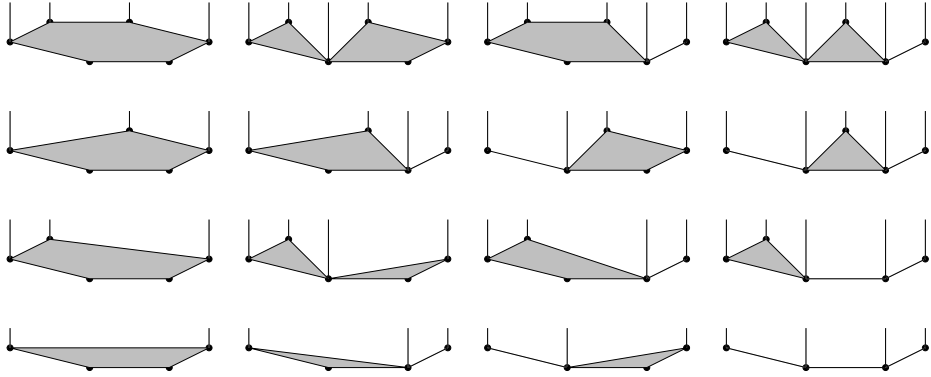}
\caption{All roofs of sub-near-edges involved in ${\overline p}_{E_{a}}$}
\label{E5a}
\end{figure}  

Figure~\ref{E5a} 
contains all roofs of the four possible 
sub-near-edges of $E_a$ obtained by removing any subset of points in
$\{P_1,P_3\}$. Their contributions to ${\overline p}_{E_a}$ are given by
$$\begin{array}{|l||r|r|r|r|}
\hline
\hbox{sub-near-edge}& &&&\\
\hline\hline
\displaystyle (P_0,P_1,P_2,P_3,P_4,P_5)&14p_3s^5&2p_4s^5&5p_4s^5&p_5s^5\\
\hline\displaystyle (P_0,P_2,P_3,P_4,P_5)&5p_2s^4&2p_3s^4&2p_3s^4&p_4s^4\\
\hline\displaystyle (P_0,P_1,P_2,P_4,P_5)&5p_2s^4&p_3s^4&2p_3s^4&p_4s^4\\
\hline\displaystyle (P_0,P_2,P_4,P_5)&2p_1s^3&p_2s^3&p_2s^3&p_3s^3
\\ \hline\end{array}$$
They sum up to the complete near-edge polynomial of $E_a$:
$${\overline p}_{E_a}=(14p_3+7p_4+p_5)s^5+(10p_2+7p_3+2p_4)s^4+
(2p_1+2p_2+p_3)s^3.$$

\subsection{Main result}

We can now state and prove the central theorem.

\begin{thm}\label{nearlyconvexcase}
  Given $l\geq 3$ near-edges  $E_1$, $E_2$... $E_l$, 
  the number of maximal triangulations of the nearly convex polygon 
  $P(E_1,\ldots,E_l)$ is given by
  $$ \tau_{\mathrm{max}}(P(E_1,\ldots,E_l)) = \big\langle \prod_{i=1}^l p_{E_i}(t),t^2G_C(t) \big\rangle_t$$
  and its triangulation polynomial is defined by
  $$ \pt(P(E_1,\ldots,E_l)) = \sum_{k}\tau_k(P(E_1,\ldots,E_l))s^k = 
      \big\langle \prod_{i=1}^l \overline{p}_{E_i}(t),t^2G_C(t) \big\rangle_t,$$  
  using the notations of 
  formulae~\ref{gencatalan}, \ref{maxnearedgepol} and~\ref{compnearedgepol} above.
\end{thm}  

\begin{cor} 
  The number of maximal triangulations and the
  triangulation polynomial of a nearly convex polygon $P(E_1,\ldots,E_l)$
  does not depend on the cyclic order of the near-edges $E_i$. 
\end{cor}

\noindent{\bf Proof of Theorem~\ref{nearlyconvexcase}.} 
Fix a triangulation $\T$ of $P=P(E_1,E_2,\ldots,E_l)$,
and a near-edge $E=E_i$. We have seen in Lemma~\ref{lempartroof} and Lemma~\ref{lemtriangunion} that 
we can associate to $E$ a partial roof $R=R(\T,E)$ and a triangulation of the sheltered
region $S(E,R)$ by the edge-type triangles of $E$.
We denote by $E'$ the subset of all points of $E$ occuring in $\T$.
The set $E'$ is a sub-near-edge since $\T$ triangulates $\conv(P)$ and thus 
involves all extremal points.
A crucial remark is that $E'\subset S(E,R)$ by Lemma~\ref{keylemma}.
Therefore, $R$ is a roof for $E'$. Moreover, by definition of $E'$,
the triangulation on $S(E,R)=S(E',R)$ is maximal with respect to $E'$. 

\begin{figure}[h]
\epsfxsize=\textwidth
\epsfbox{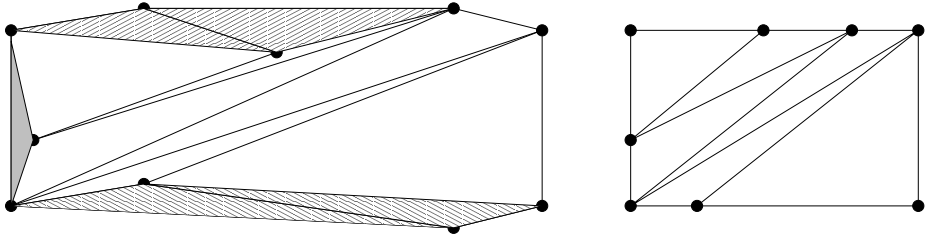}
\caption{Triangulation of the interior region $I$.}
\label{edgetriangle}
\end{figure}

Writing $R_i=R(\T,E_i)$, we obtain a triangulation of each  $S(E_i,R_i)$.
We get also a triangulation of the complement $I=\conv(P)\setminus \cup_i S(E_i,R_i)$,
whose boundary is the union of the roofs $R_i$. 
This triangulation is not arbitrary: each vertex is in a roof $R_i$ and 
a triangle in $I$ is never of edge-type and has thus not all three vertices in the same roof. 
The triangulation induced on $I$ yields hence
a triangulation of the weakly convex polygon $P(\len(R_1),\ldots,\len(R_l))$ as can be seen  
on Figure~\ref{edgetriangle}. 

We associate thus to each triangulation of $P$ a couple $((E'_i,R_i,\T_i)_i,\T_0)$, 
with the following specifications:
for each $1\leq i\leq l$,  $E'_i$ is a sub-near-edge of $E_i$, $R_i$ is a roof of $E'_i$ and $\T_i$ is 
a triangulation of  $S(E'_i,R_i)=S(E_i,R_i)$, maximal relatively to $E'i$; $\T_0$ is a maximal triangulation
of $P(\len(R_1),\ldots,\len(R_l))$.
This correspondence is clearly one-to-one: it is easy to reconstruct
a triangulation of $P$, given such data, and there is only one 
possible reconstruction.  

We consider first the case of maximal triangulations. They satisfy 
$E'_i=E_i$ for all $i$ and this property characterizes maximal triangulations. 
Theorem~\ref{weaklyconvexmaxcase} shows that the number of corresponding triangulations
equals 
$$ \left(\prod_{i=1}^l\tmax(E_i,R_i)\right)\tmax(\len(R_1),\ldots,\len(R_l)) $$
$$   = \big\langle 
       \prod_{i=1}^l\tmax(E_i,R_i)p_{\len(R_i)}(t) , t^2G_C(t)
     \big\rangle_t  $$
for each choice of a family of roofs $R_i$. 
Summation over all possible choices of $R_i$ and inversion of sum and product give the result.

For the triangulation polynomial, according to the result above,
each choice of a family of sub-near-edges gives the contribution
$$\big\langle \prod_{i=1}^l p_{E'_i}(t),t^2G_C(t) \big\rangle_t\cdot s^{\sum k_i}  ,$$
where the $k_i$ are the weights of the sub-near-edges $E'i$. Summation over all choices of sub-near-edges
and inversion of sum and product achieve the proof.
\hfill$\Box$

\noindent{\bf Example.} Consider the nearly convex polygon $P(E_a,E_b,E_c)$ 
isotopic to the perturbation depicted on Figure~\ref{exple1}. Its near-edges
$$\begin{array}{l}
\displaystyle
E_a=\left(\begin{array}{rrrrrr}0&1&2&3&4&5\\0&1&-1&1&-1&0\end{array}\right)\\
\displaystyle 
E_b=\left(\begin{array}{rrrrr}0&1&2&3&4\\0&1&-1&1&0\end{array}\right)\\
\displaystyle
E_c=\left(\begin{array}{rrrrrr}0&1&2&3&4&5\\0&2&1&-1&1&0\end{array}\right)
\end{array}$$
have complete polynomials
$$\begin{array}{lcl}
  \displaystyle
  {\overline p}_{E_a}&=&(14p_3+7p_4+p_5)s^5+(10p_2+7p_3+2p_4)s^4\\
  \displaystyle &&\quad +(2p_1+2p_2+p_3)s^3\ ,\\
  \displaystyle
  {\overline p}_{E_b}&=&(5p_3+p_4)s^4+2(2p_2+p_3)s^3+(p_1+p_2)s^2\\
  \displaystyle
  {\overline p}_{E_c}&=&
  (10p_3+7p_4+2p_5)s^5+(3p_2+13p_3+4p_4)s^4\\
  \displaystyle &&\quad +3(2p_2+p_3)s^3+(p_1+p_2)s^2,
\end{array}$$
and the triangulation polynomial
$\pt(P(E_a,E_b,E_c))$ of $P(E_a,E_b,E_c)$ is given by
$$\begin{array}{cl}
  \displaystyle &\langle {\overline p}_{E_a}(t){\overline p}_{E_b}(t)
  {\overline p}_{E_c}(t),t^2G_C(t)\rangle_t\\
  \displaystyle =&
  194939\,{s}^{14}+338669\,{s}^{13}+263615\,{s}^{12}+119944\,{s}^{11}\\
  \displaystyle 
  &\quad +34773\,{s}^{10}+6522\,{s}^{9}+748\,{s}^{8}+42\,{s}^{7}\ .
\end{array}$$


\subsection{Arithmetics of near-edges}

\subsubsection{Factorization}

A near-edge $E=(P_0,\dots,P_n)$ {\it factorizes} into near-edges
$E_1,E_2$ if there exists a lower extremal vertex
$P_k\in V^-(E)$ such that the two near-edges defined by
$$E_1=(P_0,P_1,\dots,P_{k-1},P_k),\ E_2=(P_k,P_{k+1},\dots,P_{n-1},P_n)$$
have the property that all points of $E\setminus E_i$ lie strictly above every line defined 
by two distinct points of $E_i$ for $i=1,2$. We write $E=E_1\cdot E_2$
if the near-edge $E$ factorizes with first factor $E_1$ and second 
factor $E_2$. A near-edge is {\it prime} if it has no non-trivial
factorization. It is easy to show that every near-edge 
has a unique factorization into prime near-edges.

\begin{prop} 
  Given a factorization $E=E_1\cdot E_2$ of a near-edge $E$ we have
  $$p_E=p_{E_1}\ p_{E_2}\ \ \mathrm{and}\ \ {\overline p}_E={\overline p}_{E_1}\ {\overline p}_{E_2}\ .$$
\end{prop}

\noindent{\bf Proof.} Since the 
nearly convex polygons $P(1^{k}, E)$ (with the same notation as before: $k$ successive
edges of weight $1$) and $P(1^{k}, E_1, E_2)$
are isotopic for all $k=2,3,\dots$ we have 
$$\big\langle t^k\ {\overline p}_E,t^2G_C(t)\big\rangle_t
 =\big\langle t^k\ {\overline p}_{E_1}\ {\overline p}_{E_2},t^2G_C(t)\big\rangle_t
.$$
This implies the result
since $\det(\big(C_{i+j+k}\big)_{0\leq i,j\leq n})>0$
for all $k\geq 0$ and $n\geq 1$ (this follows for instance easily from 
Exercice 6.26.b in \cite{St}).\hfill $\Box$

\metarem{Où mettre cela exactement ? Avec l'exemple numérique, mais où?}
%

\metarem{Pour mettre l'exemple numérique, il faudra d'abord mieux expliquer la philosophie}

\subsubsection{Polynomials for small near-edges}

This subsection describes all $1-$, $2-$ and $3-$near-edges up to isotopy and gives 
their polynomials. We will use the following definition:
a near-edge is {\it generic} if its underlying set of points
is a generic configuration of ${\mathbf R}^2$, i.e. if three distinct 
points of $E$ are never collinear.

We will also use the following obvious fact.
If two $n-$near-edges $E=(P_0,\dots,P_n)$ and $E'$
are {\it vertical mirrors}, \textit{i.e.} if $E'=(\overline P_n,\overline P_{n-1},\dots,\overline P_1,\overline P_0)$
where $\overline P_i=\left(\begin{array}{c}-x_i\\y_i\end{array}\right)$
is the Euclidean reflection of 
$P_i=\left(\begin{array}{c}x_i\\y_i\end{array}\right)$ with respect 
to the vertical line $x=0$,  
then $\overline p_E=\overline p_E'\ $.

\paragraph{$1-$near-edges} The unique $1-$near-edge can be
represented by $E_{1}=\left(\begin{array}{cc}0&1\\0&0\end{array}\right)$.
It is generic and prime and
has complete polynomial ${\overline p}_{E_{1}}=p_1s=s\ t$.

\paragraph{$2-$near-edges} There are two generic $2-$near-edges,
represented by
$$E_{2,1}=\left(\begin{array}{rrr}0&1&2\\0&1&0\end{array}\right),\quad
E_{2,2}=\left(\begin{array}{rrr}0&1&2\\0&-1&0\end{array}\right)\ .$$
$E_{2,1}$ is prime while $E_{2,2}=E_{1}\cdot E_{1}$. 
They have complete polynomials
$${\overline p}_{E_{2,1}}=p_2s^2+p_1s,\ {\overline p}_{E_{2,2}}=
(p_2+p_1)s^2={\overline p}_{E_{1}}^2\ .$$

Moreover, there is also a unique non-generic $2-$near-edge represented 
for instance by 
$$E_{2,3}=\left(\begin{array}{rrr}0&1&2\\0&0&0\end{array}\right)$$
with complete polynomial given by
$${\overline p}_{E_{2,3}}=p_2s^2+p_1s={\overline p}_{E_{2,1}}\ .$$

\paragraph{$3-$near-edges} There are eight
generic $3-$near-edges represented by
$$\begin{array}{ll}
\displaystyle 
E_{3,1}=\left(\begin{array}{rrrr}0&1&2&3\\0&1&3&0\end{array}\right)\quad &
E_{3,2}=\left(\begin{array}{rrrr}0&1&2&3\\0&1&1&0\end{array}\right)\ \\
\displaystyle 
E_{3,3}=\left(\begin{array}{rrrr}0&1&2&3\\0&3&1&0\end{array}\right)\ & 
E_{3,4}=\left(\begin{array}{rrrr}0&1&2&3\\0&1&-1&0\end{array}\right)\quad \\
\displaystyle 
E_{3,5}=\left(\begin{array}{rrrr}0&1&2&3\\0&-1&1&0\end{array}\right)\ &
E_{3,6}=\left(\begin{array}{rrrr}0&1&2&3\\0&-3&-1&0\end{array}\right)\quad \\
\displaystyle
E_{3,7}=\left(\begin{array}{rrrr}0&1&2&3\\0&-1&-1&0\end{array}\right)\quad &
E_{3,8}=\left(\begin{array}{rrrr}0&1&2&3\\0&-1&-3&0\end{array}\right)\ 
\quad .\\
\end{array} 
$$

The first five are prime. The last three have factorizations 
$$E_{3,6}=E_{1}E_{2,1}\ ,E_{3,7}=E_{1}^3,\ E_{3,8}=E_{2,1}E_{1}\ .$$
The pairs $\{E_{3,1},E_{3,3}\},\ \{E_{3,4},E_{3,5}\},\ \{E_{3,6},E_{3,8}\}$
are vertical mirrors. The prime near-edges have complete polynomials
$$\begin{array}{l}
\displaystyle {\overline p}_{E_{3,1}}={\overline p}_{E_{3,3}}=
(p_2+p_3)s^3+2p_2s^2+p_1s\ ,\\
\displaystyle {\overline p}_{E_{3,2}}=2p_3s^3+2p_2s^2+p_1s\ ,\\ 
\displaystyle {\overline p}_{E_{3,4}}={\overline p}_{E_{3,5}}=
(2p_2+p_3)s^3+p_1^2s^2\ .
\end{array}$$

There are moreover nine more $3-$near-edges which are not generic. They 
are represented for instance by

$$\begin{array}{ll}
\displaystyle 
E_{3,9}=\left(\begin{array}{rrrr}0&1&2&3\\0&1&2&0\end{array}\right)\quad &
E_{3,10}=\left(\begin{array}{rrrr}0&1&2&3\\0&-1&-2&0\end{array}\right)\ \\
\displaystyle 
E_{3,11}=\left(\begin{array}{rrrr}0&1&2&3\\0&0&1&0\end{array}\right)\ & 
E_{3,12}=\left(\begin{array}{rrrr}0&1&2&3\\0&0&-1&0\end{array}\right)\quad \\
\displaystyle 
E_{3,13}=\left(\begin{array}{rrrr}0&1&2&3\\0&1&0&0\end{array}\right)\ &
E_{3,14}=\left(\begin{array}{rrrr}0&1&2&3\\0&-1&0&0\end{array}\right)\quad \\
\displaystyle
E_{3,15}=\left(\begin{array}{rrrr}0&1&2&3\\0&2&1&0\end{array}\right)\quad &
E_{3,16}=\left(\begin{array}{rrrr}0&1&2&3\\0&-2&-1&0\end{array}\right)\ \quad \\
\displaystyle
E_{3,17}=\left(\begin{array}{rrrr}0&1&2&3\\0&0&0&0\end{array}\right)\quad .
\end{array} 
$$

The following near-edges factorize: 
$$E_{3,10}=E_{2,3}\ E_{1},\qquad E_{3,16}=E_{1}\ E_{2,3}$$

The remaining near-edges are prime and have complete polynomials
$$\begin{array}{lclcl}
\displaystyle 
{\overline p}_{E_{3,9}}&=&{\overline p}_{E_{3,15}}&=&p_3s^3+2p_2s^2+p_1s,\\
\displaystyle
{\overline p}_{E_{3,11}}&=&{\overline p}_{E_{3,13}}&=&
(p_2+p_3)s^3+2p_2s^2+p_1s,\\
\displaystyle 
{\overline p}_{E_{3,12}}&=&{\overline p}_{E_{3,14}}&=&(p_2+p_3)s^3+p_2s^2,\\
\displaystyle 
{\overline p}_{E_{3,17}}&=&{\overline p}_3&=&p_3s^3+2p_2s^2+p_1s\ .
\end{array}$$

\metarem{Doit-on mettre l'énumeration des near-edges generiques?}

\section{Remarks and questions}\label{raq}

\subsection{Choice of the triangulation polynomial}
 
One can also consider the triangulation
polynomial defined by
$$\sum_{k_0,k_1}\tau_{k_0,k_1}({\mathcal C})s_0^{k_0}s_1^{k_1}$$
counting the number of triangulations using $k_0$ vertices
and $k_1$ edges. The number $k_2$ of triangles can then be recovered 
using the Euler characteristic $k_0-k_1+k_2=1$ of a compact, simply 
connected triangulated polygonal region in ${\mathbf R}^2$.
This more general polynomial yields the same information as the
complete polynomial considered above except if the boundary
$\partial(\conv(\C))$ contains points of 
$\C$ which are not extremal. Most of the results
and algorithms can easily be modified in order to deal with this
more general polynomial. For clarity and concision we described here the
simpler version defined above.

\subsection{General configurations and nearly convex polygons}

Remark that every generic configuration is isotopic to
a nearly convex polygon.
Indeed, every extremal point $Q$ of a generic configuration $\mathcal C$
yields a realization of $\mathcal C$ as a nearly convex polygon
with two trivial near-edges (each consisting of $Q$ and of a
neighbouring extremal point) and
a near-edge defined by $\mathcal C\setminus \{Q\}$, which is unique
up to isotopy. 

However, the framework of near-edges is not interesting for a
general generic configuration. It speeds up computations only in the
case where the configuration has a "non-trivial" factorization
into near-edges.

\metarem{Voir avec Roland la comparaison entre itotopie et isomophisme}
%
%

\subsection{Remarks on effectiveness}

Near-edge polynomials
are in general difficult to compute.
We will present a few algorithms dealing with them in a further paper.
One of these algorithms is a slightly more sophisticated version of
an algorithm by Kaibel and Ziegler described in \cite{KZ}
and yields also a general purpose 
algorithm (unfortunately of exponential complexity), for computing
arbitrary triangulation polynomials. 
This algorithm, based on a transfer matrix, is fairly simple and it would
be interesting to compare its performance with existing 
algorithms, like for instance the algorithm of Aichholzer described
in \cite{A1}. 

The next subsection describes a family of near-edges for which
the computation of near-edge polynomials is much easier and can be achieved
by an algorithm of polynomial time-complexity.
A detailled description of the algorithm will be given in our planned future paper.

\subsection{Convex near-edges}

A near-edge $E=\{P_0,\dots,P_n\}$ 
is {\it convex} if $P_0,\dots,P_n$ are extremal points of $\conv(E)$.
Otherwise stated, the points of a convex near-edge are the vertices
of a convex polygon with $(n+1)$ edges.  There are thus
exactly $2^{n-1}$ non-isotopic convex $n-$near-edges.

A convex near-edge can be represented by a sequence of points
$$P_0=\left(\begin{array}{c}0\\ 0\end{array}\right),\dots,
P_i=\left(\begin{array}{c}i\\ \epsilon_i\ i(n-i)\end{array}\right),\dots,
P_n=\left(\begin{array}{c}n\\ 0\end{array}\right)$$
where $\epsilon_1,\dots,\epsilon_{n-1}\in\{\pm 1\}$. There are
$2^{n-1}$ equivalence classes of convex $n-$near-edges,
encoded by $n-$tuples $(\epsilon_1,\dots,\epsilon_{n-1})$ in $\{\pm 1\}^{n-1}$, The 
convex near-edge with $\epsilon_i=-1$, for all $i$, has the factorization
$E_{1}^n$. All others are prime.

A future paper will describe an algorithm having
polynomial time and memory requirements for  
computing maximal and complete edge-polynomials of convex near-edges.
It provides an efficient method for counting triangulations of nearly convex polygons
involving only convex near-edges. 
Completing each convex $n-1-$near-edge with two trivial near-edges, we get
an exponentially large
class of configurations for which the problem of counting triangulations
can be solved in polynomial time.

\subsection{Convex near-edges related to the Legendre symbol}

We used the Legendre symbol to produce data for
testing our algorithm. 
The surprising results lead to formulate the conjecture below. 

Given an odd prime $p$, the Legendre symbol, denoted by
$\left(\frac{x}{p}\right)\in \{\pm 1\}$ for $1\leq x\leq p-1$ 
defines a non-trivial homomorphism between the multiplicative groups 
$\left({\mathbf Z}/p{\mathbf Z}\right)^*$ and $\{\pm 1\}$. 
It can be computed using quadratic reciprocity or the equality
$$\left(\frac{x}{p}\right)\equiv x^{(p-1)/2}\pmod p\ .$$
We consider two convex $(p+1)-$near-edges $E_p^+,E_p^-$
associated to the sequences
$$\left(\frac{1}{p}\right),\left(\frac{2}{p}\right),\dots,\left(\frac{p-1}{p}\right)
   \ \ \mathrm{and}\ \ 
  -\left(\frac{1}{p}\right),-\left(\frac{2}{p}\right),\dots,-\left(\frac{p-1}{p}\right)$$ 
of (negated) Legendre symbols. 
For $p\equiv 3\pmod 4$ the identity
$\left(\frac{x}{p}\right)=-\left(\frac{-x}{p}\right)$ implies that
$E_p^+$ and $E_p^-$ have identical (complete) triangulation 
polynomials. 

Computation of the maximal triangulation polynomials
$p_{E_p^+}$ and $p_{E_p^-}$ for all odd primes $p<200$ suggests:

\begin{conj} Using the notations of formulae~\ref{gencatalan} and~\ref{scalprod}, we have
$$\langle P(t),t^2G_C(t)\rangle_t
    \equiv
  \left(\frac{-1}{p}\right)\pmod p=\left\{\begin{array}{cl}1&p\equiv 1\pmod 4\\
-1&p\equiv 3\pmod 4\end{array}\right.$$
for $P$ a polynomial of the form $t\ p_{E_p^+}^2,\ t\ p_{E_p^-}^2$
or $t\ p_{E_p^+}\ p_{E_p^-}$.
\end{conj}


\vspace{5mm}\noindent
{\small 
Roland Bacher and Fr\'ed\'eric Mouton, Institut Fourier\\ 
UMR 5582, Laboratoire de Math\'ematiques, BP 74 \\
F-38402 SAINT-MARTIN-D'H\`ERES CEDEX (FRANCE)  \\
Roland.Bacher@ujf-grenoble.fr\\
Frederic.Mouton@ujf-grenoble.fr
}

\end{document}